\documentclass[reqno,12pt]{amsart}
\usepackage{amsthm}
\usepackage{amsmath, mathrsfs}
\usepackage{color}
\usepackage{amssymb}
\usepackage{hyperref}
\usepackage{enumerate}
\usepackage{latexsym,array}
\usepackage{amsfonts}
\usepackage{shadow}

\newtheorem{Pa}{Paper}[section]
\newtheorem{Tm}[Pa]{{\bf Theorem}}
\newtheorem{La}[Pa]{{\bf Lemma}}

\newtheorem{Cy}[Pa]{{\bf Corollary}}

\newtheorem{Rk}[Pa]{{\bf Remark}}
\newtheorem{Pn}[Pa]{{\bf Proposition}}

\newtheorem{Ex}[Pa]{{\bf Example}}

\newcommand{\w}{\omega}
\def\e{\epsilon_N}

\author[D. Alpay]{Daniel Alpay}
\address{(DA) Department of Mathematics\newline
BenGurion University of the Negev\newline Beer-Sheva 84105
Israel} \email{dany@math.bgu.ac.il}

\author[P. Jorgensen]{Palle Jorgensen}
\address{(PJ)
Department of Mathematics\newline 14 MLH \newline The University
of Iowa\newline Iowa City, IA 52242-1419 USA}
\email{palle-jorgensen@uiowa.edu}

\author[I. Lewkowicz]{Izchak Lewkowicz}
\address{(IL) Department of Electrical Engineering
\newline
Ben Gurion University of the Negev \newline P.O.B. 653,
\newline
Be'er Sheva 84105, \newline ISRAEL} \email{izchak@ee.bgu.ac.il}

\title[realizations of infinite products and wavelet filters]{
Realizations of infinite products,
Ruelle operators and wavelet filters} \oddsidemargin 0.2in
\def\R{\mathbb R}

\def\(s){\mathscr S(\R\times\R)}



\thanks{{\sl Acknowledgments:} D. Alpay thanks the Earl Katz family for
endowing the chair which supported his research. The authors
thank the US-Israel Binational Science Foundation (BSF) Grant number 2010117.}

\begin{document}
\maketitle \tableofcontents
\parindent 0cm
\begin{abstract}
Using the system theory notion of state-space realization of
matrix-valued rational functions, we describe the Ruelle
operator associated with wavelet filters.
The resulting realization of infinite products of rational
functions have the following four features: 1) It is
defined in an infinite-dimensional complex domain.  2) Starting
with a realization of a single rational matrix-function $M$, we
show that a resulting infinite product realization obtained from
$M$ takes the form of an (infinite-dimensional) Toeplitz operator
with the symbol that is a reflection of the initial realization for
$M$.  3)  Starting with a subclass of rational matrix functions,
including scalar-valued ones corresponding to low-pass wavelet
filters, we obtain the corresponding infinite products that
realize the Fourier transforms of generators of $\mathbf
L_2(\mathbb R)$ wavelets.   4) We use both the realizations for
$M$ and the corresponding infinite product to obtain a matrix
representation of the Ruelle-transfer operators used in wavelet
theory.  By ``matrix representation'' we refer to the slanted (and
sparse) matrix which realizes the Ruelle-transfer operator under
consideration. 
\end{abstract}

\parindent 0cm

\section{Introduction}
\setcounter{equation}{0}
Among many applications of rational matrix-valued functions are
their use as filters in signal processing, and in the
construction of classes of wavelets. In the latter case, the
matrix function to be considered is made up of a prescribed
system of scalar valued functions of a single complex variable.
If $N$ is a scaling number for the wavelet under consideration,
then there are associated systems of $N$ scalar-valued functions
representing each of the corresponding $N$ frequency bands. Each
such system produces a matrix-valued function. This particular
approach to wavelet filters was considered in
\cite{BrJo02b,BrJo02a,ajlm1,ajl1paper,GNS, MR3050315}. The function
corresponding to low-pass yields a father function for a wavelet
when certain technical assumptions are imposed. Here we consider
instead the matrix-valued approach: it has the advantage that it
allows one to treat the combination of individual bands in a
single analysis. However the issues involving infinite products
in the matrix-valued case are more subtle, and we address them
below. For example, to understand the infinite product formed
from a rational matrix-valued function of a single complex
variable, one must introduce an infinite number of complex
variables. We show that, under suitable assumptions, the infinite
product-function in turn then also has a realization as a
function of one variable. While our motivation derived initially
from the study of wavelet filters, we note that there is a host
of other applications of infinite products of rational matrix
functions. Indeed, the framework for our consideration of
infinite products goes beyond that of wavelet filters. We shall
consider these more general settings in the last section of our
paper. The latter non-wavelet applications are derived from the
theory of systems. Indeed the theory of realization of systems is
also a key tool in our analysis of
infinite products.\\

In \cite[Section 4]{ajl1paper}, we characterized wavelet
filters as functions of the form
\begin{equation}
\label{waveletfilter} M(z)=QU(z^N)\Delta(z)V
\end{equation}
where 
\[
V={\scriptstyle\frac{1}{\sqrt{N}}}\left(\e^{-\ell j}\right)_{\ell,j=0,\ldots,
N-1}\quad\quad\quad\quad
\e:=e^{i\frac{2\pi}{N}},
\]
is (up to scaling) the usual discrete Fourier transform matrix,
\[
\Delta(z):=\left(\begin{smallmatrix}
1&      &    &      &         \\
 &z^{-1}&    &      &         \\
 &           &\ddots&         \\
 &           &      &z^{-(N-1)}\end{smallmatrix}\right),
\]
$U$ is a rational $(N\times N)$-valued function which takes unitary
values on the unit circle, with no poles outside the closed unit
disk, and and $Q$ is an arbitrary (constant) unitary matrix.
One can explicitly write \eqref{waveletfilter} as
\begin{equation}\label{ajlmmmm}
M(z)={\scriptstyle\frac{1}{\sqrt{N}}}
\begin{pmatrix}m_0(z)    &m_0(\e z)  &\cdots&m_0(\e^{N-1}z)\\
               m_1(z)    &m_1(\e z)  &\cdots&m_1(\e^{N-1}z)\\
                \vdots   &           &      &             \\
               m_{N-1}(z)&m_{N-1}(\e z)&\cdots& m_{N-1}(\e^{N-1}
               z)\end{pmatrix}.
\end{equation}
Note that $M(z)$ in \eqref{waveletfilter} is unitary on $\mathbb{T}$.\\

An earlier relevant result on wavelet filter \eqref{waveletfilter},
\eqref{ajlmmmm} appeared in \cite{GNS}. In \cite{AJL3} and
\cite{AJL4}, we recently further explored rectangular rational
functions which are (co)isometric on the unit circle (with poles
anywhere, but $\mathbb{T}$).\\

Following \eqref{ajlmmmm} one can write,
\begin{equation}
\begin{pmatrix}m_0(z)\\m_1( z)\\ \vdots\\
m_{N-1}(z)\end{pmatrix}=QU(z^N)\begin{pmatrix}1\\z^{-1}\\
\vdots\\z^{-(N-1)}\end{pmatrix}. \label{qwert1111}
\end{equation}
In the sequel, by choosing in \eqref{waveletfilter}
\begin{equation}\label{eq:Q}
Q=\left(U(1)\Delta(1)V\right)^*=\left(U(1)V\right)^*
\end{equation}
we shall normalize the filters so that in \eqref{waveletfilter}
\[
M(1)=I_N~.
\]
This normalization in particular forces that the upper left
entry of $M$ to satisfy
\[
m_0(1)=1.
\]
This last condition is crucial to consider infinite products.
For $m_0(z)$ in \eqref{ajlmmmm} we set:
\begin{equation}\label{eq:m0}
m(z):=m_{0}(z).
\end{equation}
The wavelet father function $\varphi(w)$ is
given by its Fourier transform
\begin{equation}
\label{infiniteprod}
\widehat{\varphi}(w)=\prod_{k=1}^\infty
m(e^{\frac{2\pi iw}{N^k}}).
\end{equation}
For details, see e.g. \cite{BrJo02a} and \cite{BrJo02b}.\\

It should be pointed out that in some engineering circles 
functions of the form of \eqref{waveletfilter} are referred to
as (multi-resolution) {\em Filter Bank}, whose applications
transcend wavelets, see e.g. \cite{TV}, \cite{Va}
and even go beyond signal processing, see e.g. \cite{LPV}.\\

Here we limit our
discussion to discrete time systems, where the variables,
input, output, and state, are time series, that is, functions on
$\mathbb Z$. A linear time-invariant system in this model will
then be specified by a transfer matrix  $M(z)$, also called a
transfer function. It is a rational matrix-valued function of a
single complex variable. Moreover the complex variable $z$ is dual
to time, and so it represents frequency. If $M(z)$ has no pole
at infinity (following engineering literature), a corresponding
state space realization is any quadruple of matrices $A, B, C$, and
$D$ of appropriate size,
such that
\[
M(z)=D+C(zI-A)^{-1}B
\]
holds.

We assume that $M$ in \eqref{waveletfilter}, \eqref{ajlmmmm} is a
matrix-valued rational function analytic at infinity,
while for $m(z)$, its upper left entry, we introduce
a state space realization
\begin{equation}\label{eq:m1}
m(z)=D+C(zI-A)^{-1}B
\end{equation}
where we can assume that the realization is minimal (that is, the
size of $A$ is the smallest possible one), and that in
particular $A$ has no spectrum on the unit circle since $M$ is
analytic on the unit circle. Our use of the term realization
conforms to its common use in the theory of systems from the study
of dynamical systems and filters in engineering, and pioneered by
Kalman and others; see \cite{MR569473, MR0255260, MR1640001}.
Aspects of realization of filter bank as in \eqref{waveletfilter},
\eqref{ajlmmmm} were already addressed in \cite{GNS}, \cite{TV}
and \cite{Va}.\\

The paper consists of five sections besides the introduction,
and its outline is as follows. In Section \ref{sec:FiniteProd}
we introduce state-space realization formulas of finite
products of rational functions, each of a different variable.
Infinite products are considered in Section \ref{sec4}. As we
will see in that section an important role is played by the
Toeplitz operator with the related symbol equal to
\begin{equation}
A+zB(I-zD)^{-1}C. \label{symboltoeplitz}
\end{equation}
In Section \ref{sec5}, we compute the Markov parameters associated
with $|m(z)|^2$ in terms of the given realization of $m$. In Section
\ref{sec6}, we study the Ruelle operator and connections with
rational wavelet filters are studied in Section \ref{sec7}. 

\section{Finite products}
\setcounter{equation}{0} \label{sec:FiniteProd} 

As is well known, see e.g  \cite{bgk1}, \cite[Section 6.4]{MR569473}, 
\cite{MR0255260}, \cite[Section 6.5]{MR1640001}, every
$(p\times q)$-valued rational function $R(z)$ 
analytic at infinity can be written as
\begin{equation}
R(z)=D+C(zI-A)^{-1}B, \label{real1}
\end{equation}
for matrices $A,B,C$ and $D$ of appropriate sizes. Equation
\eqref{real1} is called a realization of $R$ and we shall
sometimes use the abbreviated form
\begin{equation}\label{eq:realiz}
\left(\begin{array}{c|c}A&B\\
\hline
C&D\end{array}\right).
\end{equation}
In general, the realization is highly non unique.
When the dimension of the state space, (i.e. the upper
left block $A$ is say $d\times d$), is {\em minimal}, the
realization is unique up to a similarity matrix, meaning that
the only freedom in the choice of the realization is
\begin{equation}
\label{eq:similarity}
\begin{pmatrix}A&B\\ C&D\end{pmatrix}\,\,\mapsto\,\,
\begin{pmatrix}T&0\\ 0&I_p\end{pmatrix}
\begin{pmatrix}A&B\\ C&D\end{pmatrix}
\begin{pmatrix}T^{-1}&0\\ 0&I_q\end{pmatrix},
\end{equation}
where $T\in\mathbb C^{d\times d}$ is an arbitrary invertible
matrix and where $D$ is assumed to belong to $\mathbb C^{p\times
q}$. An important formula for the realization of the product is
given in the next lemma:
\begin{La}
Let $R_1$ and $R_2$ be two matrix-valued rational functions
analytic at infinity, and with realizations
\[
R_1(z)=D_1+C_1(zI_{n_1}-A_1)^{-1}B_1\quad and\quad
R_2(z)=D_2+C_2(zI_{n_2}-A_2)^{-1}B_2.
\]
Assume that the product $R_1R_2$ makes sense. Then a realization
of $R_1(z)R_2(z)$ is compactly given by:
\begin{equation}\label{Alma-Marceau}
\left(\begin{array}{cc|c}
A_1&B_1C_2&B_1D_2\\~
&A_2&B_2\\
\hline
C_1&D_1C_2&D_1D_2
\end{array}\right)
\end{equation}
\end{La}
An important tool in our argument is the counterpart of
\eqref{Alma-Marceau}  when each function depends on a 
{\it different
variable}. See Lemma \ref{Varenne}.\\

Factorization of rational matrix-valued functions of one variable
is classical. In contrast, factorization theory of rational
functions of several complex variables $z_1,~\ldots~,~z_u$ is
not well developed. However, here we consider matrix-valued
rational functions of the form
\begin{equation}\label{qwertyu}
M(z)=M_1(z_1)M_2(z_2)\cdots M_u(z_u)
\end{equation}
where $M_1,~\ldots~,~M_u$ are matrix-valued rational functions of
appropriate sizes and analytic at infinity.\\

For future reference we mention the following result, whose
proof is a direct verification, and will be omitted.
\begin{La}
Let $\mathcal H_1$ and $\mathcal H_2$ be two Hilbert spaces and
let $a\, :\,\mathcal H_1\,\,\rightarrow\,\,
\mathcal H_1$, $b\, :\,\mathcal H_2\,\,\rightarrow\,\,
\mathcal H_2$, and  $c\, :\,\mathcal H_2\,\,\rightarrow\,\,
\mathcal H_1$ be bounded linear operators, with $a$ and $b$ invertible. Then
\[
\begin{pmatrix}a&-c\\0&b\end{pmatrix}^{-1}=
\begin{pmatrix}a^{-1}&a^{-1}cb^{-1}\\0&b^{-1}\end{pmatrix}.
\]
\label{la3-1}
\end{La}

The following very simple lemma is the key to the formulas we develop:
\begin{La}
\label{Varenne}
\begin{equation}
\label{eq:formula1}
\begin{split}
\left(D_1+C_1(z_1I_{n_1}-A_1)^{-1}B_1\right)
\left(D_2+C_2(z_2I_{n_2}-A_2)^{-1}B_2\right)=\\
&\hspace{-4cm}=D+C(\Lambda(z)-A)^{-1}B
\end{split}
\end{equation}
where, 
\begin{equation}\label{Lambda}
\Lambda(z):=\begin{pmatrix}z_1I_{n_1}&~\\~&z_2I_{n_2}\end{pmatrix}
\end{equation}
and the realization array \eqref{eq:realiz} takes the form,
\begin{equation}
\left(\begin{array}{cc|c}
A_1&B_1C_2&B_1D_2\\~
&A_2&B_2\\
\hline
C_1&D_1C_2&D_1D_2
\end{array}\right)
\end{equation}
\end{La}
\begin{proof}
We first note that by Lemma \ref{la3-1} we have:
\[
\begin{split}
(\Lambda(z)-A)^{-1}&=
\begin{pmatrix}z_1I_{n_1}-A_1&-B_1C_2\\
0&z_2I_{n_2}-A_2\end{pmatrix}^{-1}\\
&=
\begin{pmatrix}(z_1I_{n_1}-A_1)^{-1}&
(z_1I_{n_1}-A_1)^{-1}B_1C_2(z_2I_{n_2}-A_2)^{-1}\\
0&(z_2I_{n_2}-z_2A_2)^{-1}\end{pmatrix}.
\end{split}
\]
Therefore we have
\[
\begin{split}
D+C(\Lambda(z)-A)^{-1}B&=D_1D_2+\\
&\hspace{-5cm}+\begin{pmatrix}C_1&D_1C_2\end{pmatrix}
\begin{pmatrix}(z_1I_{n_1}-z_1A_1)^{-1}&
(z_1I_{n_1}-A_1)^{-1}B_1C_2(z_2I_{n_2}-A_2)^{-1}\\
0&(z_2I_{n_2}-A_2)^{-1}\end{pmatrix}\begin{pmatrix}B_1D_2\\
B_2\end{pmatrix}
\end{split}
\]
which is exactly the left side of \eqref{eq:formula1}.
\end{proof}
This formula can now be iterated to obtain a realization for a
product \eqref{qwertyu}. With
\[
M_j(z_j)=D_j+C_j(z_jI_{n_j}-A_j)^{-1}B_j,\quad j=1,\ldots u,
\]
and $M(z)=M_1(z_1)M_2(z_2)\cdots M_u(z_u)$ we have
\[
M(z)=D+C(\Lambda(z)-A)^{-1}B,
\]
where
\[
\Lambda(z):=\left(\begin{smallmatrix}
z_1I_{n_1}&          &      &         \\
          &z_2I_{n_2}&      &         \\
          &          &\ddots&         \\
          &          &      &z_uI_{n_u}
\end{smallmatrix}\right)
\]
where $M(\infty)=\lim\limits_{z\rightarrow\infty}M(z)$
and the realization array \eqref{eq:realiz} takes the form,
\begin{equation}\label{newu}
\left(\begin{array}{cccccc|c}
A_1&B_1C_2  &  B_1D_2C_3&\cdots &B_1D_2\cdots D_{u-2}C_{u-1}&
B_1D_2\cdots D_{u-1}C_u&B_1D_2\cdots D_u\\
~  & A_2    &  B_2C_3   & \cdots &B_2D_3\cdots D_{u-2}C_{u-1}&
B_2D_3\cdots D_{u-1}C_u&B_2D_3\cdots D_u\\
~&~&~&\ddots&~&~&\vdots \\
~&~&~&~&A_{u-1}&B_{u-1}C_u&B_{u-1}D_u \\
~  &~     &~        &~&~      &A_u       & B_u\\
\hline
C_1&D_1C_2&D_1D_2C_3&\cdots&D_1\cdots D_{u-2}C_{u-1}&
D_1\cdots D_{u-1}C_u)&D_1\cdots D_u
\end{array}\right)
\end{equation}
We note that the case where all the functions vanish at infinity, i.e.
where \mbox{$\lim\limits_{z\rightarrow\infty}M(z)=0$,}
leads to very simple formulas, which we gather in the
following lemma.

\begin{La}\label{LemmaD=0}
It holds that
\[
C_1(z_1I_{n_1}-A_1)^{-1}B_1C_2(z_2I_{n_2}-A_2)^{-1}B_2
\cdots C_u(z_uI_{n_u}-A_u)^{-1}B_u=C(\Lambda(z)-A)^{-1}B,
\]
where $\Lambda(z)$ is as in \eqref{Lambda}
and the realization array \eqref{eq:realiz} takes the form,
\[
\left(\begin{array}{cccccc|c}
A_1&B_1C_2  &0&\cdots &~&0&0\\
~  & A_2    &  B_2C_3   &0& \cdots &0&0\\
~&~&~&\ddots&~&~&\vdots \\
~&~&~&~&A_{u-1}&B_{u-1}C_u&0\\
~  &~     &~        &~&~      &A_u       & B_u\\
\hline
C_1&0&~&\cdots&~&0&0
\end{array}\right).
\]
\end{La}
The significance of the realizations in \eqref{newu} and in
Lemma \ref{LemmaD=0}, goes beyond the scope of this work.
In the sequel, we actually exploit a special case of it. See
also Remark \ref{remark}

\section{Infinite products}
\setcounter{equation}{0} \label{sec4} While the framework of
realizations is typically formulated for finite matrices, (as we
point out below) a number of the results make sense for infinite
matrices, hence for linear operators in Hilbert space. A case in
point is the realizations we obtain now for our infinite products.
We now wish to let $u\rightarrow\infty$ in \eqref{newu} when all
the functions $R_j$ coincide:
\[
R_1(z)=R_2(z)=\cdots=M(z),
\]
where $M(z)$ is a matrix-valued rational function, analytic at
infinity, with realization
\begin{equation}\label{eq:m}
M(z)=D+C(zI-A)^{-1}B.
\end{equation}
We assume that $1\geq\| M(z)\|$ for all $z\in\mathbb T$.

\begin{Tm}Given a square $M(z)$ in \eqref{eq:m}.\mbox{}\\
$(i)$
Assume that
\begin{equation}
\label{ddd}
\lim\limits_{k\rightarrow\infty}\|D^k\|^{1/k}<1.
\end{equation}
Then, the operators
\begin{eqnarray}
\label{eq:newop1}
\mathscr A&=&\begin{pmatrix} A&BC&BDC&BD^2C&\cdots\\
0&A&BC&BDC&\cdots\\
0&0&A&BC&\cdots\\
& & & & \\
 & & & &
\end{pmatrix},\quad \ell_2(\mathbb N)\otimes \mathbb
C^m\,\,\Longrightarrow\,\, \ell_2(\mathbb N)\otimes \mathbb C^m\\
\label{eq:newop2}
\mathscr B&=&\begin{pmatrix}\vdots\\ BD^2\\ BD\\
B\end{pmatrix},\hspace{2.3cm}\quad \mathbb C^m
\longrightarrow\,\, \ell_2(\mathbb N)\otimes\mathbb C^m,\\
\mathscr C&=&\begin{pmatrix}
C&DC&D^2C&\cdots\end{pmatrix},\hspace{0.3cm} \quad \ell_2(\mathbb
N)\otimes \mathbb C^m \longrightarrow\,\,\mathbb C^m,
\label{eq:newop3}
\end{eqnarray}
are bounded.\smallskip

$(ii)$ $\mathscr A$ in \eqref{eq:newop1} is the block-Toeplitz
operator with symbol
\[
A+zB(I-zD)^{-1}C.
\]
\end{Tm}

\begin{proof} We have
\[
A+zB(I-zD)^{-1}C=A+zBC+z^2BDC+z^3BD^2C+\cdots,
\]
and hence the function $\phi(z)=A+zB(I-zD)^{-1}C$ is the symbol of
the block Toeplitz operator \eqref{eq:newop1}. We note that, in
view of \eqref{ddd}
\begin{equation}
\label{eq:toepliz}
\|\mathscr A\|=
\|A+zB(I-zD)^{-1}C\|_\infty<\infty,
\end{equation}
and so the block Toepliz operator $\mathscr A$ is bounded. We use
the fact that a block Toeplitz operator with symbol $\phi(z)$ has
norm $\|\phi\|_\infty$; see \cite{ggk1}.
\end{proof}

\begin{Rk}{\rm
The same result holds {\it mutatis mutandis} when instead of
\eqref{eq:m}, the alternative
realization $M(z)=D+zC(I-zA)^{-1}B$ is chosen.}
\end{Rk}

\begin{Rk}{\rm  Recall that
\[
\mathscr B=\mathcal{O}(D, B)\quad{\rm  and}\quad
\mathscr C=\mathcal{C}(D, C)
\]
in \eqref{eq:newop2} and \eqref{eq:newop3} are the observability
and controllability operators (see e.g. \cite[Section 6.2]{MR569473},
\cite{MR1640001}) 
\label{Rk4.2}
}
\end{Rk}

\begin{Tm}
Assume that $M(z)$ in \eqref{eq:m} has no singularity at the point
$z=1$ and that $M(1)=I$, and let $(z_k)_{k\in\mathbb N_0}$ be a
sequence of complex numbers
which are not poles of $M$ and such that
\begin{equation}
\label{St-Jacques}
\sum_{k=0}^\infty |1-z_k|<\infty.
\end{equation}
Then it holds that
\begin{equation}
\label{Sully-Morland-Ligne-7} \prod_{k=1}^\infty M(z_k)=\mathscr
C(\Lambda(z)-\mathscr A)^{-1}\mathscr B
\end{equation}
where $\mathscr A,\mathscr B$ and $\mathscr C$ are defined by
\eqref{eq:newop1}-\eqref{eq:newop3}, and where
\[
\Lambda(z)=\left(\begin{matrix}
z_1I_n&~     &~     \\
~     &z_2I_n&      \\
~     &~     &\ddots
\end{matrix}\right).
\]
\end{Tm}

\begin{proof} Since the realization of $M$ is assumed minimal, $1$
is not in the spectrum of $A$ and
we have $M(1)=D+C(I_n-A)^{-1}B$. Thus
\[
\begin{split}
M(z)-M(1)&=D+C(zI_n-A)^{-1}B-D-C(I_n-A)^{-1}B\\
&=C(zI_n-A)^{-1}B-C(I_n-A)^{-1}B\\
&=(1-z)C(zI-A)^{-1}(I-A)^{-1}B.
\end{split}
\]
Furthermore, \eqref{St-Jacques} implies in particular that
$\lim_{k\rightarrow\infty}z_k=1$. Let
\[
K_0=\max_{z\in V}\left\{\|C(zI_n-A)^{-1}(I_n-A)^{-1}B\|\right\},
\]
where $V$ is a closed neighborhood of $1$ in which $m$ has no
pole, and let
\[
K=\max_{\text{where the $z_u\not\in V$}}\left\{K_0,\|
C(z_uI_n-A)^{-1}(I_n-A)^{-1}B\|\right\}
\]
Then we have
\[
\| M(z_k)-I\|\le K\cdot|1-z_k|.
\]
Therefore the series $\sum_{k=1}^\infty \| I-M(z_k)\|$ and hence
the product $\prod_{k=1}^\infty M(z_k)$ are convergent. The
equality \eqref{Sully-Morland-Ligne-7} is now easy to verify.
\end{proof}

\begin{Cy}
Assume that $M(1)=I$, and let $(\theta_k)$ be a sequence of
numbers on the real line such that
\[
\sum_{k=0}^\infty |\theta_k|<\infty.
\]
Then the infinite product
\[
\prod_{k=0}^\infty M(e^{it\theta_k}),\quad t\in\mathbb R,
\]
converges for all real $t$. \label{les-thibault}
\end{Cy}

\begin{proof}
Since $|e^{i\theta}-1|\le |\theta|$ for $\theta$ real, we have
\[
\| M(e^{it\theta_k})-I\|\le K_1\cdot|\theta_k|
\]
where now we can take
\[
K_1=\max_{\theta\in[0,2\pi]}\|C(e^{i\theta}I_n-A)^{-1}(I_n-A)^{-1}B\|,
\]
and hence the result.
\end{proof}

\begin{Cy}
In the notation and hypothesis of Corollary \ref{les-thibault}, the product
\[
\prod_{k=0}^\infty M(e^{i\frac{2\pi w}{N^k}})
\]
converges for every $w\in\mathbb R$.
\end{Cy}

\section{Markov parameters}
\setcounter{equation}{0} \label{sec5}
Let $M$ be analytic in the exterior and on the boundary of the unit disk, with
minimal realization $M(z)=D+C(zI-A)^{-1}B$. In particular
$\sigma(A)\subset\mathbb D$. Let $h_0+\sum_{k=1}^\infty
\frac{h_k}{z^k}$ be the Laurent expansion at infinity of $M$. The
coefficients $h_0,h_1,\ldots$ are called the Markov parameters of
$M$, see e.g. \cite[Subsections 5.1.2, 6.2.1]{MR569473},
\cite[Section 6.5]{MR1640001}. They are given by $h_0=D$ and
\begin{equation}
h_k=CA^{k-1}B,\quad k=1,2,\ldots \label{Tolbiac}
\end{equation}
We extend the sequence $h_k$ by
\[
h_u=0,\quad u<0.
\]
Since the spectral radius of $A$ is strictly less than
$1$, we can set
\begin{equation}
\label{gamma}
\Gamma=\sum_{u=0}^\infty A^{*u}C^*CA^u.
\end{equation}
Note that $\Gamma$ is called the observability Gramian, and is the
unique solution of the Stein equation
\[
\Gamma-A^*\Gamma A=C^*C.
\]
We set
\begin{equation}
\label{M} Y=D^*C+B^*\Gamma A.
\end{equation}

In view of the next result we recall that a rational function $r$
with no poles on the unit circle belongs to the Wiener
algebra of the disk, that is, can be written as
\[
r(z)=\sum_{n\in\mathbb Z}z^nr_n,
\]
where $\sum_{n\in\mathbb Z}|r_n|<\infty$. See for instance
\cite[Corollary 3.2]{gk1}.

\begin{Tm}
Let $(c_n)_{n\in\mathbb Z}$ be defined by
\begin{equation}
\label{m2} \| M(z)\|^2=\sum_{n\in\mathbb Z}c_nz^n,\quad z\in\mathbb
T.
\end{equation}
Then
\[
c_n=\sum_{j\in\mathbb Z}h_j^*h_{j+n},\quad n\in\mathbb Z. \]
Therefore,
\begin{equation}
c_n=\begin{cases}\,\, B^*A^{*(-n-1)}Y^*,\,\,\,\hspace{2mm}
n<0,\\
\,\, D^*D+B^*\Gamma B,\quad\hspace{0.1cm} n=0,\\
\,\, YA^{n-1}B,\,\, \hspace{1.2cm}n>0.
\end{cases}
\label{formulasRuelle}
\end{equation}
\label{thmruelle}
\end{Tm}

\begin{proof}
For $n=0$ we have
\[
\begin{split}
c_0&=D^*D+\sum_{k=1}^\infty B^*A^{*(k-1)}C^*CA^{k-1}B\\
&=D^*D+ B^*\left(\sum_{k=1}^\infty A^{*(k-1)}C^*CA^{k-1}\right)B\\
&=D^*D+B^*\Gamma B,
\end{split}
\]
where $\Gamma$ is the Gramian matrix from \eqref{gamma}.\\

We now assume $n>0$. Then,
\[
\begin{split}
c_n&=h_0^*h_n+\sum_{k=1}^\infty h_k^*h_{k+n}\\
&=D^*CA^{n-1}B+\sum_{k=1}^\infty B^*A^{k-1}C^*CA^{k+n-1}B\\
&=D^*CA^{n-1}B+B^*\Gamma A^nB.
\end{split}
\]
Finally, for $n<0$, we have:
\[
\begin{split}
c_n&=\sum_{k=-n}^\infty h_k^*h_{k+n}\\
&=\sum_{u=0}^\infty h_{u-n}^*h_u\\
&=B^*A^{*(-n-1)}C^*D+\sum_{u=1}^\infty
B^*A^{*(u-n-1)}C^*CA^{u-1}B\\
&=B^*A^{*(-n-1)}C^*D+B^*A^{*(-n)}\Gamma B.
\end{split}
\]
\end{proof}

\begin{Cy}
Let $d$ be the size of $A$ (that is, $A\in\mathbb C^{d\times d}$).
Let $Y$ be defined by \eqref{M} and
\[
\mathcal C(A,B)=\begin{pmatrix}B&AB&A^{2}B&\cdots&A^{d-1}B\end{pmatrix}.
\]
Then
\begin{equation}
c_n=YA^{n-1}B,\quad n=1,\ldots 
\label{CH}
\end{equation}
and
\begin{equation}
\begin{pmatrix}c_1&c_2&\cdots &c_{d}\end{pmatrix}= Y\cdot
\mathcal C(A,B).
\end{equation}
\end{Cy}

\begin{Tm}
There exists complex numbers $a_0,\ldots, a_{d-1}$ such that
\begin{equation}
\label{recursion} a_0c_1+a_1c_2+\cdots +a_{d-1}c_{d}+c_{d+1}=0,
\end{equation}
and more generally, for any $p\ge 1$,
\begin{equation}
\label{recursion2} a_0c_p+a_1c_{p+1}+\cdots +a_{d-1}c_{d+p-1}+c_{d+p}=0.
\end{equation}
\end{Tm}
\begin{proof}
By the Cayley-Hamilton theorem there exists numbers
$a_0,a_1,\ldots , a_{d-1}$ such that
\[
a_0+a_1A+a_2A^{2}+\cdots +a_{d-1}A^{d-1}+A^d=0.
\]
It then follows from \eqref{CH} that we have \eqref{recursion},
and more generally \eqref{recursion2}.
\end{proof}

Formulas \eqref{formulasRuelle} take a simpler form in a number of
cases, which we mention as remarks:

\begin{Rk}
{\rm If $A$ is nilpotent (and then $m$ is a polynomial in $1/z$)
we have
\[
c_{d+1}=c_{d+2}=\cdots =0.
\]
}
\end{Rk}

\begin{Rk}
{\rm Assuming that $D=0$ implies that in \eqref{M} $Y=D^*C$
and \eqref{formulasRuelle} becomes,
\begin{equation}
c_n=\begin{cases}\,\, B^*A^{*(-n)}\Gamma B,\,\,
\,n\le 0,\\
\,\,B^*\Gamma A^{n}B,\hspace{.65cm}\,\,\,n>0.
\end{cases}
\label{formulasRuelle2}
\end{equation}
}
\end{Rk}

\begin{Rk}{\rm
We recall that the observability Gramian $\Gamma$, see
\eqref{gamma}, is invertible
if and only if the pair $(C,A)$ is observable, meaning that
\[
\cap_{u=1}^{d-1}\ker A^u=\left\{0\right\}
\]
Then one assume $\Gamma=I_d$ by taking $T=\Gamma^{1/2}$ as
similarity matrix in \eqref{eq:similarity}. When furthermore $D=0$
we then have
\begin{equation}
c_n=\begin{cases}\,\, B^*A^{*(-n)} B,\,\,
\,n\le 0,\\
\,\,B^* A^{n}B,\hspace{.65cm}\,\,\,n>0.
\end{cases}
\label{formulasRuelle23}
\end{equation}
 }
\end{Rk}

\begin{Pn}
Assume $M$ is rational. Then the coefficients $c_n$ satisfy the
estimates of the form
\begin{equation}
\label{ruelle} |c_{k}|\le Ce^{-\alpha|k|},\quad k\in\mathbb Z.
\end{equation}
for every $\alpha>\rho(A)$, where $\rho(A)$ is the spectral
radius of $A$.
\end{Pn}

\section{Ruelle operator}
\setcounter{equation}{0} \label{sec6} 
Using $m(z)$ from \eqref{eq:m0}, \eqref{eq:m1}, the associated
Ruelle operator (or transfer operator) is defined by (see
\cite[\S 3.2]{BrJo02a})
\begin{equation}
(Rf)(z)={\scriptstyle\frac{1}{N}}\sum_{\substack{w\in\mathbb T
\\w^N=z}}|m(w)|^2f(w).
\end{equation}
See \cite[p. 156]{BrJo02a}. In terms of the coefficients
\eqref{m2} in Theorem \ref{thmruelle} it is the operator between
appropriate subspaces of $\ell_2(\mathbb Z)$ and with matrix
representation
\[
r_{\ell,j}={\scriptstyle\frac{1}{N}}c_{N\ell-j},
\]
that is
\begin{equation}
\label{rfl}
(Rf)_\ell={\scriptstyle\frac{1}{N}}\sum_{k\in\mathbb Z}
c_{N\ell-k}f_k.
\end{equation}
\begin{Ex}
{\rm
When $N=2$, $D=0$ and $\Gamma=I_d$ the matrix representation of
the Ruelle operator is 
\[
\frac{1}{2}\begin{pmatrix}
\cdots &      &        &           &          &  &        &\cdots \\
\cdots &B^*B  & B^*AB  & B^*A^{*2}B&B^*A^{*3}B&B^*A^{*4}B&&\cdots\\
\cdots &B^*A^2B  & B^*AB  &
\fbox{$B^*B$}&B^*A^{*}B&B^*A^{*2}B&B^*A^{*3}B&\cdots\\
\cdots &B^*A^4B  & B^*A^3B  &
B^*A^2B&B^*AB&B^*B&B^*A^{*}B&\cdots\\
\cdots &      &        &           &          & &         & \cdots\\
\end{pmatrix}
\]
where the box denotes the $(0,0)$ element (see \eqref{formulasRuelle23}).
}
\end{Ex}

Let
\[ \mathscr E_r=\left\{(f_n)_{n\in\mathbb Z}\,;\,
\|f\|_{r,1}\stackrel{\rm def.}{=}
\sum_{n\in\mathbb Z}e^{r|n|}|f_n|<\infty\right\}
\]
and
\[ \mathscr E_r^{(2)}=\left\{(f_n)_{n\in\mathbb Z}\,;\,
\|f\|_{r,2}^2\stackrel{\rm def.}{=}\sum_{n\in\mathbb
Z}e^{r|n|}|f_n|^2<\infty\right\}
\]
See \cite{MR1041203}, \cite[p. 158]{BrJo02a} for these last
spaces. We note that an element of $\mathscr E_r$ satisfies an
estimate of the form
\begin{equation}
\label{ineqfn} 
|f_n|\le Ke^{-r|n|},\quad n\in\mathbb Z,
\end{equation}
for some $K>0$.

\begin{Tm} Assume that the coefficients $c_n$ satisfy the estimate
\eqref{ruelle}. Then for every choice of $\beta$ and
$\beta^\prime$ such that
\begin{equation}
\label{abbprime} \alpha<\beta\quad and\quad \beta^\prime<N\alpha,
\end{equation}
the Ruelle operator is continuous from $\mathscr E_\beta$ into
$\mathscr E_{\beta^\prime}$ and from $\mathscr E_\beta^{(2)}$ into
$\mathscr E_{\beta^\prime}^{(2)}$.
\end{Tm}
\begin{proof}  Let $f=(f_n)_{n\in\mathbb Z}$ be an element of
$\mathscr E_\beta$, with $\beta$ as in \eqref{abbprime}. The
$\ell$ component of the vector $Rf$ is given by \eqref{rfl},
\[
(Rf)_\ell={\scriptstyle\frac{1}{N}}
\sum_{k\in\mathbb Z}c_{N\ell-k}f_k
\]
and, using \eqref{ruelle} and \eqref{ineqfn} can be estimated as:
\[
\begin{split}
N|(Rf)_\ell|&\le\sum_{k\in\mathbb Z}|c_{N\ell-k}|\cdot|f_k|\\
&\le \sum_{k\in\mathbb Z}Ce^{-\alpha|N\ell-k|}\cdot K e^{-\beta
|k|}\\
&= \sum_{k\in\mathbb Z}CKe^{-\alpha(|N\ell-k|+|k|)}\cdot
e^{(\alpha-\beta)
|k|}\\
&\le CKe^{-\alpha|N\ell|}\sum_{k\in\mathbb Z}e^{(\alpha-\beta)
|k|}.
\end{split}
\]
Hence for any $\beta^\prime<N\alpha$,
\[
\sum_{\ell\in\mathbb Z}|(Rf)_\ell)|\cdot e^{\beta^\prime |\ell|}<\infty,
\]
and so $Rf\in\mathscr E_{\beta^\prime}$.\\

Let now $\epsilon>0$ and $f,g\in\mathscr E_\beta$ and such that
$\|f-g\|_{\beta,1}<\epsilon$. Then in a way similar to the
above argument we write
\[
\begin{split}
N\|f-g\|_{\beta^\prime,1}&=\sum_{\ell\in\mathbb Z}\left(
\sum_{k\in\mathbb Z}|c_{N\ell-k}|\cdot|f_k-g_k|\right)\\
&\le \sum_{\ell\in\mathbb Z}\left(\sum_{k\in\mathbb Z}
Ce^{-\alpha|N\ell-k|}\cdot  e^{-\alpha |k|}\underbrace{
e^{(\alpha-\beta^\prime)|k|}}_{\text{is $<1$}}
e^{\beta^\prime |k|}|f_k-g_k|\right)\\
&\le \sum_{\ell\in\mathbb Z}\left(\sum_{k\in\mathbb Z}
Ce^{-\alpha|N\ell|}\cdot e^{\beta^\prime |k|}|f_k-g_k|\right)\\
&=C\left(\sum_{\ell\in\mathbb Z}e^{-\alpha|N\ell|}\right)
\cdot\|f-g\|_{\beta^\prime,1}.
\end{split}
\]
and continuity of $R$ follows.\\

The case of the spaces $\mathscr E^{(2)}_{\beta}$ and $\mathscr
E^{(2)}_{\beta^\prime}$ is proved in the same way.
\end{proof}

For a result related to the following theorem in the non rational
case, see \cite[p. 158]{BrJo02a}.\\

The first item in the next theorem is taken from
\cite[p.156-159]{BrJo02a}, \cite{Dau95}.
\begin{Tm}Let $m(z)$ be as in \eqref{eq:m0}, \eqref{eq:m1}.\mbox{}\\
$(1)$ The Ruelle operator has finite trace, and its trace is given
by the formula
\[
{\rm Tr}\, R\,={\scriptstyle\frac{1}{N}}\sum_{k\in\mathbb
Z}c_{(N-1)k}={\scriptstyle\frac{1}{N}}\sum_{\substack{w\in\mathbb T
\\ w^N=1}}|m(w)|^2.
\]
$(2)$ In the rational case, we have
\[
{\rm Tr}\, R\,=\,{\scriptstyle\frac{1}{N}}
\left(DD^*+B^*\Gamma B+Y(I-A)^{-1}B^*+B(I-A^*)^{-1}Y^*\right).
\]
\end{Tm}

\begin{proof}
The second item is a direct consequence of formulas
\eqref{formulasRuelle} .
\end{proof}

\section{Wavelets and rational filters}
\setcounter{equation}{0} \label{sec7} In this section, we show that
starting from a rational wavelet filter, the infinite product
\eqref{infiniteprod} is indeed in $\mathbf L_2(\mathbb R,dx)$. To
that purpose it is enough to prove that
\[
R1\le 1
\]
for the corresponding Ruelle operator, where, by definition of $R$,
\[
(R1)(z)={\scriptstyle\frac{1}{N}}
\sum_{\substack{w\in\mathbb T\\w^N=z}}|m(w)|^2,
\]
with $m$ from \eqref{eq:m0}, \eqref{eq:m1}.

Let $M$ be a rational wavelet filter, that is a function of the
form \eqref{waveletfilter}. Recall that its first column is given
by \eqref{qwert1111}.

\begin{Pn} It holds that
\begin{equation}
R1=1
\end{equation}
\end{Pn}

\begin{proof}
From \eqref{qwert1111} and \eqref{eq:Q} we have that
\[
m(z)=\sum_{j=1}^Nz^{1-j}[W(z^N)]_{1j},
\]
where $[W]_{jk}$ denotes the $jk$ entry in a matrix $W$, with
\[
W(z^N):=(U(1)V)^*U(z^N).
\]
Thus,
for $z\in\mathbb T$,
\[
\begin{split}
{\scriptstyle\frac{1}{N}}\sum_{\substack{\omega\in\mathbb T\\
{\omega}^N=z}}|m(\omega)|^2&={\scriptstyle\frac{1}{N}}
\sum_{\substack{\omega\in\mathbb T\\{\omega}^N=z}}\left|
\sum_{j=1}^N{\omega}^{1-j}\left[W({\omega}^N)\right]_{1j}
\right|^2\\
&={\scriptstyle\frac{1}{N}}\sum_{\substack{\omega\in\mathbb T\\
{\omega}^N=z}}\sum_{j,k=1}^N{\omega}^{-j+k}
\left[W(z)\right]_{1j}\overline{\left[W(z)\right]_{1k}}
\quad\quad(\mbox{\rm note that ${\omega}^N=z$})\\
&=
\sum_{j,k=1}^N\delta_{jk}
\left[W(z)\right]_{1j}\overline{\left[W(z)\right]_{1k}}
\quad\quad(\mbox{\rm
since
${\scriptstyle\frac{1}{N}}\sum\limits_{\substack{\omega\in\mathbb T\\
{\omega}^N=z}}{\omega}^{-j+k}=\delta_{jk}$})\\
&=\sum_{j=1}^N
\left[W(z)\right]_{1j}\overline{\left[W(z)\right]_{1j}}
=1\quad\quad(\mbox{\rm since $W$ takes
unitary values on $\mathbb T$}).
\end{split}
\]
\end{proof}

\begin{Tm}
Let $m$ be the upper left entry of a wavelet filter,
see \eqref{eq:m0}, \eqref{eq:m1}. Then the
infinite product \eqref{infiniteprod} converges to an $\mathbf
L_2(\mathbb R)$ function.
\end{Tm}

\begin{proof} We  proceed in a number of steps.\\

STEP 1: {\sl The infinite product \eqref{infiniteprod} converges
pointwise.}\smallskip

This follows from Corollary \ref{les-thibault} with
$\theta_k=\frac{2\pi w}{N^k}$ since $M(1)=I_N$ and hence
$m(1)=1$.\\

We now set
\begin{equation}
\label{fk-1}
f_k(w)=1_{[-\frac{N^k}{2},\frac{N^k}{2}]}(w)\prod_{\ell=0}^k
m\left(e^{i\frac{{2\pi}w}{N^\ell}}\right),\quad k=0,1,\ldots
\end{equation}
The key to the proof is to establish the identity
\begin{equation}
\label{keyformula} \int_{-\frac{N^k}{2}}^{\frac{N^k}{2}}|
f_k(w)|^2dw=1.
\end{equation}

STEP 2: {\sl \eqref{keyformula} holds:}\smallskip

Indeed,
\[
\begin{split}
\int_{-\frac{N^k}{2}}^{\frac{N^k}{2}}| f_k(w)|^2dw&=
\int_{-\frac{N^k}{2}}^{\frac{N^k}{2}}\left(\prod_{\ell=0}^k|m
\left(e^{i\frac{{2\pi}w}{N^\ell}}\right)|^2\right)dw\\
&=\int_{-\frac{N^{k-1}}{2}}^{\frac{N^{k-1}}{2}} \left(
\prod_{\ell=0}^{k-1}
|m\left(e^{i\frac{{2\pi}w}{N^\ell}}\right)|^2\right)
\underbrace{\left(\sum_{p=0}^{N-1}|m\left(e^{{2\pi}i
\frac{w+pN^{k-1}}{N^k}}\right)|^2\right)}_{
\mbox{\text{This is $R1$}}} dw\\
&= \int_{-\frac{N^{k-1}}{2}}^{\frac{N^{k-1}}{2}} \left(
\prod_{\ell=0}^{k-1}
|m\left(e^{i\frac{{2\pi}w}{N^\ell}}\right)|^2\right)dw\quad
(\text{since $R1=1$})\\
&=\int_{-\frac{N^{k-1}}{2}}^{\frac{N^{k-1}}{2}}
|f_{k-1}(w)|^2dw\quad (\text{by definition of $f_{k-1}$; see \eqref{fk-1}}).
\end{split}
\]
So we have
\[
\int_{-\frac{N^k}{2}}^{\frac{N^k}{2}}| f_k(w)|^2dw=
\int_{-\frac{N^{k-1}}{2}}^{\frac{N^{k-1}}{2}} |f_{k-1}(w)|^2dw.
\]
Iterating this equality we get to
\[
\begin{split}
\int_{-\frac{N^k}{2}}^{\frac{N^k}{2}}| f_k(w)|^2dw=
\int_{-\frac{N^{k-1}}{2}}^{\frac{N^{k-1}}{2}} |f_{k-1}(w)|^2dw=\cdots\\
&\hspace{-4cm}\cdots=
\int_{-\frac{1}{2}}^{\frac{1}{2}} |f_{0}(w)|^2dw=
\int_{-\frac{1}{2}}^{\frac{1}{2}} |m(e^{2\pi iw})|^2dw<\infty
\end{split}
\]

STEP 3: {\sl The pointwise limit function $f(w)=
\prod_{\ell=1}^\infty m\left(e^{2\pi i\frac{w}{N^\ell}}\right)$
belongs to $\mathbf L_2(\mathbb R)$.}\smallskip

Indeed, by Fatou's lemma,
\[
\begin{split}
\int_{\mathbb R}|f(w)|^2dw&=
\int_{\mathbb R}\lim_{k\rightarrow\infty}|f_k(w)|^2dw\\
&=\int_{\mathbb R}\liminf_{k\rightarrow\infty}|f_k(w)|^2dw\\
&\le\liminf_{k\rightarrow\infty}\int_{\mathbb
R}|f_k(w)|^2dw<\infty
\end{split}
\]
in view of the previous step.
\end{proof}

\begin{Rk}{\rm
The preceding arguments hold still in the case $R1\le 1$. This
covers important cases of rational filters for which the function
$U$ in  \eqref{waveletfilter} is only contractive as opposed to
unitary.}
\end{Rk}

\begin{Rk}\label{remark}{\rm
In  future work we shall use the multivariable state space techniques
introduced in Section \ref{sec:FiniteProd} to extend our results to
the multivariable case, including the case of an infinite number of
variables, and to study connections with infinite dimensional analysis
(see \cite{MR2387368} for the latter).}
\end{Rk}

{\bf Acknowledgments:} It is a pleasure to thank the referees for a very careful reading of the manuscript.
\bibliographystyle{plain}
\def\cprime{$'$} \def\lfhook#1{\setbox0=\hbox{#1}{\ooalign{\hidewidth
  \lower1.5ex\hbox{'}\hidewidth\crcr\unhbox0}}} \def\cprime{$'$}
  \def\cfgrv#1{\ifmmode\setbox7\hbox{$\accent"5E#1$}\else
  \setbox7\hbox{\accent"5E#1}\penalty 10000\relax\fi\raise 1\ht7
  \hbox{\lower1.05ex\hbox to 1\wd7{\hss\accent"12\hss}}\penalty 10000
  \hskip-1\wd7\penalty 10000\box7} \def\cprime{$'$} \def\cprime{$'$}
  \def\cprime{$'$} \def\cprime{$'$} \def\cprime{$'$}

\end{document}